\documentclass[]{elsarticle}

\usepackage{lineno,hyperref}

\usepackage{graphicx}
\usepackage{bm}
\usepackage{multirow}
\usepackage{amsmath}
\usepackage{amsfonts}
\usepackage{amssymb}
\usepackage{amsthm}
\usepackage{secdot}
\usepackage{tabularx}
\usepackage{xfrac}
\usepackage{longtable}
\usepackage{tablefootnote}

\journal{Applied Numerical Mathematics}

\theoremstyle{plain}

\theoremstyle{definition}

\theoremstyle{remark}

\numberwithin{equation}{section}
\numberwithin{theorem}{section}
\numberwithin{remark}{section}

\DeclareFontFamily{U}{mathx}{\hyphenchar\font45}
\DeclareFontShape{U}{mathx}{m}{n}{
      <5> <6> <7> <8> <9> <10>
      <10.95> <12> <14.4> <17.28> <20.74> <24.88>
      mathx10
      }{}
\DeclareSymbolFont{mathx}{U}{mathx}{m}{n}
\DeclareFontSubstitution{U}{mathx}{m}{n}
\DeclareMathAccent{\widecheck}{0}{mathx}{"71}
\DeclareMathAccent{\wideparen}{0}{mathx}{"75}

\begin{document}

\begin{frontmatter}

\title{An Explicit Construction of Optimized Interpolation Points on the 4-Simplex}

\author{Trenton J. Gobel \corref{mycorrespondingauthor}}

\cortext[mycorrespondingauthor]{Corresponding author}
\ead{tjg5291@psu.edu}

\author{David M. Williams} 

\address{Department of Mechanical Engineering, The Pennsylvania State University, University Park, Pennsylvania 16802}

\fntext[fn1]{Distribution Statement A: Approved for public release. Distribution is unlimited.}

\begin{abstract} 
In this work, a family of symmetric interpolation points are generated on the four-dimensional simplex (i.e.~the pentatope). These points are optimized in order to minimize the Lebesgue constant. The process of generating these points closely follows that outlined by Warburton in ``An explicit construction of interpolation nodes on the simplex," Journal of Engineering Mathematics, 2006. Here, Warburton generated optimal interpolation points on the triangle and tetrahedron by formulating explicit geometric warping and blending functions, and applying these functions to equidistant nodal distributions. The locations of the resulting points were Lebesgue-optimized. In our work, we extend this procedure to four dimensions, and construct interpolation points on the pentatope up to order ten. The Lebesgue constants of our nodal sets are calculated, and are shown to outperform those of equidistant nodal distributions. 
\end{abstract}

\begin{keyword}
Interpolation \sep space-time \sep four dimensions \sep finite element methods \sep Lebesgue constant \sep symmetry
\MSC[2010] 76M10 \sep 65D05 \sep 65M12 \sep 65M60
\end{keyword}

\end{frontmatter}

\section{Introduction} \label{intro_sec}

The nodal polynomial approximation of a well-defined function is called a `nodal representation'. One may use nodal representations in order to approximate \emph{all} functions for a given finite element method (FEM), including both the solution and the geometry. Such FEM's are referred to as nodal FEM's. However, alternative options exist. In particular, one may use either nodal or modal polynomials to represent the FEM solution and geometry. With this in mind, one may consider four different types of FEM's: a) FEM's which use a nodal representation for both the solution and geometry, b) FEM's which use a modal representation for the solution and a nodal representation for the geometry, c) FEM's which use a nodal representation for the solution and a modal representation for the geometry, and d) FEM's which use a modal representation for both the solution and geometry. We note that FEM's of types c) and d) are not very practical, as it is difficult to approximate the geometry with a modal representation. In particular, one must obtain modal coefficients for the geometry by solving linear systems of equations; these systems must be customized for each element in the mesh, and the sizes of the systems grow rapidly with the polynomial order and dimension. In contrast, nodal representations bypass this issue, as they directly interpolate the geometry. Consequently, there is no need to solve a customized linear system for the nodal coefficients, for each element in the mesh. As a result, FEM's of types a) and b) are the most popular in practice. 

We can summarize the above discussion by stating that virtually all popular FEM's require nodal interpolation points in some capacity. 

In what follows, we will discuss the important characteristics of nodal and modal representations, the key differences between the two, the unique challenges associated with constructing a nodal representation, and the importance of the Lebesgue constant.

\section{Background}

Prior to discussing the specifics of the nodal representation, it is important to first provide some contextual background by describing the modal representation and its properties. 

In order to construct a modal representation for the solution, we multiply modal coefficients by a set of element-spanning polynomials—i.e.~modes. These modes are not required to be orthogonal, but are often chosen as such in order to facilitate implementation. In particular, the orthogonality of the polynomials results in a diagonal mass matrix thereby simplifying computations. The modal coefficients and polynomials are independent of any interpolation node locations, and confer the advantage that changing the polynomial order does not necessarily require changing all of the basis functions. This latter property is guaranteed to hold if the modal basis is hierarchical.

To construct a nodal representation, we must choose specific point locations in the reference element at which the values of the solution will be defined or sampled. These solution values are then multiplied by nodal basis functions (usually Lagrange polynomials) and summed in order to return the solution on the entire element. Many researchers consider nodal representations to be more intuitive than modal representations, as the nodal coefficients directly correspond to the values of the FEM solution at particular points in physical space. Conversely, modal coefficients correspond to the amounts of energy that are stored in particular modes. Unfortunately, nodal basis functions are neither hierarchical nor orthogonal, which makes their implementation/computation less straightforward than their modal counterparts. 

Naturally, this raises the question of how nodal basis functions can be generated. In one dimension, these basis functions can be calculated analytically using the well-known Lagrange formula. The extension of this approach to higher dimensions is straightforward on quadrilaterals, hexahedra, tesseracts and other hypercube elements, as the tensor products of Lagrange polynomials can be used to construct multi-dimensional nodal basis functions. However, the construction of nodal basis functions on higher-dimensional simplices is more difficult, owing in part to the nonexistence of an elegant, tensor-product-based procedure. Some explicit basis function formulae exist—such as those given by Hughes~\cite{hughes2012finite} for triangles and tetrahedra—but these functions are limited to low polynomial orders, and are difficult to generate in an automatic fashion. 

Evidently, a different approach for generating nodal basis functions is needed. Fortunately, as discussed in~\cite{Hesthaven07}, Vandermonde matrices relate modal representations to nodal representations. The modal basis in question is constructed using products of Jacobi polynomials of the desired order. Utilizing these Jacobi polynomials, we can construct an orthonormal modal basis on the reference simplicial element. Once this is obtained, we can multiply a vector of the modal basis functions by the inverse transpose of the Vandermonde matrix, and thereafter determine the values of the multi-dimensional Lagrange polynomials.

Proper behavior of the Vandermonde matrix is crucial for interpolation accuracy and stability, especially at high polynomial orders. Therefore, one often seeks to maximize the absolute value of the matrix determinant. The maximization process depends on two factors. Firstly, we must ensure the orthonormality of the modal polynomial basis. Fortunately, an equation for an orthonormal basis in four dimensions on the pentatope has already been obtained—see the formulation in~\cite{Frontin21}. 
Secondly, we must select interpolation points for the construction of the Vandermonde matrix which are optimally distributed throughout the simplex. It turns out that Legendre-Gauss-Lobatto points do an excellent job of amplifying the Vandermonde matrix determinant in one dimension. Therefore, we seek to construct generalizations of Legendre-Gauss-Lobatto points on to higher-dimensional simplices.

In addition to maximizing the determinant of the Vandermonde matrix, there are two qualities that bear consideration: symmetry and error minimization. Symmetry of interpolation points on the simplex is defined as invariance of the point locations under affine transformations of the simplex on to itself. Error minimization is defined based on the points' ability to minimize the interpolation error with respect to a particular norm. 

It is well-known that minimization of the pointwise error is equivalent to minimizing the error with respect to the $L^\infty$ norm. In turn, minimizing the $L^\infty$ norm is achieved by minimizing the Lebesgue constant for a given nodal distribution. Finally, it is important to note that minimizing the Lebesgue constant is directly connected to maximizing the Vandermonde matrix determinant. Therefore, we can focus on minimizing the Lebesgue constant, and obtain some of the numerical-stability benefits associated with maximizing the determinant without additional effort.

As an alternative strategy, we can also minimize the amount of error introduced during numerical quadrature procedures. We will refer to this strategy as quadrature-error minimization. It turns out, that for low polynomial orders, quadrature-optimized interpolation points can have relatively small Lebesgue constants. As a result, they simultaneously maintain good integration accuracy and $L^{\infty}$ accuracy~\cite{Williams13}. 

In summary, one often seeks to construct interpolation points on the simplex which satisfy the following objectives: i) the points are arranged in a symmetric fashion, ii) the points maximize the determinant of the Vandermonde matrix, and iii) the points minimize either quadrature error or $L^{\infty}$ error.  

\subsection{Literature Review}

Early efforts to construct interpolation points focused on maximizing the determinant of the Vandermonde matrix. Interpolation points constructed in this fashion are called Fekete points~\cite{Fekete23}. In the 1980's, Bos used this optimization approach to generate symmetric interpolation points on triangles up to order seven~\cite{Bos83}. Chen and Babuška furthered this concept up to order twenty on the line segment and order thirteen on the triangle~\cite{Chen95}. They also constructed symmetric interpolation points on the tetrahedron up to order nine~\cite{chen1996optimal}. Thereafter, Taylor, Wingate, and Vincent generated symmetric interpolation points up to order nineteen on the triangle~\cite{Taylor00}, using the same methodology. In~\cite{Heinrichs05}, Heinrichs used a different approach, minimizing the Lebesgue constant directly rather than maximizing the Vandermonde matrix determinant. With this approach, interpolation points up to order eighteen were generated on the triangle. Curiously, the resulting points were not symmetric.

Another optimization method was formulated based on electrostatic principles by Hesthaven~\cite{Hesthaven98}. In this work, an analogy was drawn between the electric charge distributions of Stieltjes~\cite{Stieltjes1,Stieltjes2} and the Legendre-Gauss-Lobatto points in one dimension. Hesthaven extended the electrostatic analogy into two dimensions, and constructed symmetric interpolations points on the triangle with favorable Lebesgue constants up to order twelve.  

In addition, Warburton proposed a simple optimization procedure for minimizing the Lebesgue constant based on geometric movement (warping and blending) of point locations~\cite{Warburton06}. In this procedure, the interpolation points are initially distributed on an equispaced grid within the simplex, and then allowed to move in conjunction with explicitly defined warping and blending functions. The functions are defined on the faces and edges of the simplex, and are parameterized by a single variable. In turn, this variable is optimized in order to minimize the Lebesgue constant. Using this approach, Warburton generated symmetric interpolation points on the triangle up to order fifteen, and on the tetrahedron up to order ten.

Recently, Jameson et al.~\cite{jameson2012non} and Castonguay et al.~\cite{castonguay2011application} proved the effectiveness of quadrature-optimized interpolation points for FEM's applied to non-linear advection problems. In this work, they leveraged symmetric quadrature rules which happened to have the correct number of points to represent polynomials on the line segment and triangle. Thereafter, Shunn and Ham~\cite{shunn2012symmetric}, and Williams et al.~\cite{williams2014symmetric} expanded this work, and constructed symmetric quadrature-optimized interpolation points on the triangle and tetrahedron up to orders ten and six, respectively. In~\cite{Williams13}, Williams and Jameson showed that the Lebesgue constants for these points are well-behaved for up to order five on the tetrahedron. Most recently, Williams et al.~\cite{Williams20} extended the work of~\cite{williams2014symmetric}, and generated quadrature-optimized, symmetric interpolation points on the pentatope up to order five.  

\subsection{Overview of the Paper}

The majority of previous work has focused on generating interpolation points in one, two, and three dimensions. The comparatively small amount of work that focuses on four-dimensional space has been limited to quadrature-optimized points. In this paper, we seek to generate new, Lebesgue-optimized interpolation points on 4-simplex elements (pentatopes). Our work is a natural extension of Warburton's approach~\cite{Warburton06} into four dimensions. To our knowledge, this represents the first ever such endeavor. 

In section~\ref{prelim_sec} of this paper, we discuss the preliminary concepts and formulae relevant to this method. Section~\ref{method_sec} explains the methodology behind the warping and blending functions used in the generation of the Lebesgue-optimized points. The final results of our optimization process are given in section~\ref{results_sec}. Lastly, conclusions are given in section~\ref{conclusion_sec} along with suggestions for future work. Tables of the optimized interpolation points are given in an appendix subsequent to the main text. 

\section{Preliminaries} \label{prelim_sec}

\subsection{Nodal Representations} \label{nodal_rep_sec}

A nodal solution representation in $d$-dimensional space is given by 
\begin{align}
    u^h(\bm{r}) = \sum_{i=1}^{N_p}u_i^hl_i(\bm{r}),
\end{align}
where $u^h$ is the finite element solution, $N_p$ is the number of nodes, $l_i$ are the nodal basis functions,  and $\bm{r}$ is the coordinate vector in reference space. Here, $u_i^h$ represents the value of the solution at the specific nodal locations $\bm{r}_i$ in reference space. For a given dimension $d$ and polynomial order $p$, the number of nodes is calculated via
\begin{align}
    N_p = \frac{1}{d!}\prod_{k=1}^d(p+k).
\end{align}
In one dimension, this simplifies to $N_p = p + 1$. The nodal basis functions can be calculated for a one-dimensional mesh using Lagrange polynomials
\begin{align}
	l_j(x) = \prod^{N_p}_{\stackrel{m = 0}{m \neq j}}\frac{x-x_m}{x_j-x_m} = \frac{x-x_0}{x_j-x_0}\cdot\cdot\cdot\frac{x-x_{j-1}}{x_j-x_{j-1}}\frac{x-x_{j+1}}{x_j-x_{j+1}}\cdot\cdot\cdot\frac{x-x_{N_p}}{x_j-x_{N_p}},
\end{align}
where each $x_j$ represents a specific nodal location. However, as mentioned previously, it is difficult to construct analytical formulas for nodal basis functions on simplex elements in $\mathbb{R}^d$ when $d\geq 2$. This issue will be addressed in the next section.

\subsection{Modal Representations} \label{modal_rep_sec}

The general form for modal solution representations is very similar to nodal representations
\begin{align}
    u_h(\bm{r}) = \sum_{i=1}^{N_p}\hat{u}_i^h\psi_i(\bm{r}),
\end{align}
where $\hat{u}_i^h$ are the modal coefficients, and $\psi_i(\bm{r})$ are the modal basis functions.
However, unlike nodal basis functions, there exist many analytical formulas for modal basis functions. Oftentimes, these formulas are written in terms of Jacobi polynomials. We note that the Jacobi polynomials themselves can be formulated recursively as follows
\begin{align}
    P_0^{(\alpha,\beta)}(x) &= \sqrt{2^{-\alpha-\beta-1}\frac{\Gamma(\alpha+\beta+2)}{\Gamma(\alpha+1)\Gamma(\beta+1)}}, \\[1.0ex]
    P_1^{(\alpha,\beta)}(x) &= \frac{1}{2}P_0^{(\alpha,\beta)}(x)\sqrt{\frac{\alpha+\beta+3}{(\alpha+1)(\beta+1)}}((\alpha+\beta+2)x+(\alpha-\beta)), \\[1.0ex]
    P_{n+1}^{(\alpha,\beta)}(x) &= \frac{1}{a_{n+1}}\left(xP_n^{(\alpha,\beta)}(x)-a_nP_{n-1}^{(\alpha,\beta)}(x)-b_nP_n^{(\alpha,\beta)}(x)\right),
\end{align}
where $\alpha, \beta$ and $n$ are integer parameters, $\Gamma(\cdot)$ is the well-known Gamma function, and $a_n$ and $b_n$ are defined as
\begin{align}
    a_n &= \frac{2}{2n+\alpha+\beta}\sqrt{\frac{n(n+\alpha+\beta)(n+\alpha)(n+\beta)}{(2n+\alpha+\beta-1)(2n+\alpha+\beta+1)}}, \\[1.0ex]
    b_n &= -\frac{\alpha^2-\beta^2}{(2n+\alpha+\beta)(2n+\alpha+\beta+2)}.
\end{align}
Hesthaven and Warburton~\cite{Hesthaven07} give a Jacobi polynomial-based formulation in one dimension for the orthonormal modal basis as follows
\begin{align}
    \psi_m(r) = \sqrt{\frac{2m-1}{2}}P_{m-1}^{(0,0)}(r), \qquad 0 \leq m \leq p + 1.
\end{align}
Note that here the $\alpha$ and $\beta$ parameters of the Jacobi polynomial are zero, returning the $(m-1)^{\text{th}}$ Legendre polynomial in place of the Jacobi polynomial. 

Now that we have an expression for the modal basis, we can relate it to the nodal basis. This action is completed via  the Vandermonde matrix $\mathcal{V}$
\begin{align}
    \mathcal{V}_{ij} = \psi_j(\bm{r}_i), \qquad \mathcal{V}^T\bm{l(r)} = \tilde{\boldsymbol{\psi}}\bm{(r)},
\end{align}
where $\bm{r}_i$ are the nodal point locations in reference space, (as mentioned previously).
Consequently,
\begin{align}
    \bm{l(r)} = \mathcal{V}^{-T}\tilde{\boldsymbol{\psi}}\bm{(r)}. \label{basis}
\end{align}
Therefore, in order to obtain a higher-dimensional nodal basis, we merely need to construct a higher-dimensional \emph{modal} basis, and then use Eq.~\eqref{basis}. For the sake of completeness, in what follows we provide the modal basis formulas for two, three, and four dimensions.

In two dimensions, the modal basis is given by
\begin{align}
    \psi_m(\bm{r}) = \sqrt{2}P_i^{(0,0)}(a)P_j^{(2i+1,0)}(b)(1-b)^i,
\end{align}
where $m$ is the modal index for two-dimensional polynomials, defined over $i$ and $j$ as
\begin{align}
    m = 1 + \left(\frac{2p+3}{2}\right)i + \left(-\frac{1}{2}\right)i^2 + j, \quad (i,j) \geq 0, \quad i + j\leq p,
\end{align}
and $a$ and $b$ are defined relative to the reference coordinates $\bm{r} = (r,s)$
\begin{align}
    a = 2\left(\frac{1+r}{1-s}\right)-1, \quad b = s.
\end{align}
The modal basis in three dimensions is an extension of the two-dimensional case
\begin{align}
    \psi_m(\bm{r}) = 2\sqrt{2}P_i^{(0,0)}(a)P_j^{(2i+1,0)}(b)P_k^{(2i+2j+2,0)}(c)(1-b)^i(1-c)^{i+j},
\end{align}
where again we have a modal indexing function
\begin{align}
    m &= 1 + \left(\frac{3p^2+12p+11}{6}\right)i + \left(\frac{-p-2}{2}\right)i^2 + \left(\frac{1}{6}\right)i^3 + \left(\frac{2p+3}{2}\right)j \notag\\[1.0ex] 
    &+ \left(-\frac{1}{2}\right)j^2 - ij + k, \quad (i,j,k)\geq 0, \quad i + j + k\leq p,
\end{align}
and for $\bm{r} = (r,s,t)$
\begin{align}
    a = -2\left(\frac{1+r}{s+t}\right)-1, \quad b=2\left(\frac{1+s}{1-t}\right)-1, \quad c=t.
\end{align}

Frontin et.~al.~\cite{Frontin21} give a formulation for the four-dimensional modal basis which follows the pattern set by its lower-dimensional predecessors
\begin{align}
    \psi_{m}(\bm{r}) = 8 P_i^{(0,0)}(a) P_j^{(2i+1,0)}(b)P_k^{(2i+2j+2,0)}(c)P_q^{(2i+2j+2k+3,0)}(e) \notag \\
     \times  (1-b)^i(1-c)^{i+j}(1-e)^{i+j+k},
\end{align}
and for $\bm{r} = (r,s,t,u)$
\begin{align}
    a = -2\left(\frac{1+r}{s+t+u+1}\right)-1, \quad b=-2\left(\frac{1+s}{t+u}\right)-1, \quad c=2\left(\frac{1+t}{1-u}\right)-1, \quad e = u.
\end{align}
As in the lower dimensional cases, a modal indexing formula is necessary in four dimensions. However, no expression for $m$ has been derived in previous works. 

By observation, we realize that the indexing equation is polynomial in nature with each term of the form $C_ni^{\xi} j^{\zeta} k^{\mu} q^{\nu}$, where $\xi$, $\zeta$, $\mu$, and $\nu$ are non-negative integers. We also note that each coefficient $C_n$ is either constant or a polynomial function of $p$. This suggests the following formula for $m$ in four dimensions
\begin{align}
    m &= 1 + C_1i + C_2i^2 + C_3i^3 + C_4i^4 + C_5j + C_6ij + C_7i^2j + C_8j^2 + C_9ij^2+ C_{10}j^3 \notag \\[1.0ex] 
    &+ C_{11}k + C_{12}ik + C_{13}jk +  C_{14}k^2 + C_{15}q, \notag \\[1.0ex] 
    &(i,j,k,q)\geq 0, \quad i+j+k+q\leq p. \label{m_partial} 
\end{align}
Note that every permutation of $i,j,k,$ and $q$ is not represented in the above formula. Certain terms were excluded based on the assumption that, like the lower-dimensional equations before it, there exists an upper limit on the order of a term depending on the variables appearing in that term and the dimension $d$. For example, no term containing an $i$ can have an order above $d$, $j$ above $d-1$, $k$ above $d-2$, and $q$ above $d-3$. 

To find explicit expressions for the coefficients in Eq.~\eqref{m_partial}, we must first generate the correct node index ordering using a series of nested \textit{for} loops. Once proper ordering is ascertained, we carefully choose values of $i,j,k,$ and $q$ to calculate the numerical values for each variable $C_n$ for a given $p$ as shown in Table~\ref{C_Table_num}. The MATLAB code used to generate Table~\ref{C_Table_num} is given in Appendix B of \cite{Gobel24}.
\begin{longtable}{c | c c c c c}
\caption{Four-dimensional node indexing equation coefficients for varying $p$} \label{C_Table_num} \\
   \hline
   $p$ & 4 & 5 & 6 & 7 & 8 \\
   \hline
   \hline
   $C_1$ & $\sfrac{533}{12}$ & $\sfrac{275}{4}$ & $\sfrac{1207}{12}$ & $\sfrac{1691}{12}$ & $\sfrac{763}{4}$\\
   $C_2$ & $\sfrac{-251}{24}$ & $\sfrac{-335}{24}$ & $\sfrac{-431}{24}$ & $\sfrac{-539}{24}$ & $\sfrac{-659}{24}$\\
   $C_3$ & $\sfrac{13}{12}$ & $\sfrac{5}{4}$ & $\sfrac{17}{12}$ & $\sfrac{19}{12}$ & $\sfrac{7}{4}$ \\
   $C_4$ & $\sfrac{-1}{24}$ & $\sfrac{-1}{24}$ & $\sfrac{-1}{24}$ & $\sfrac{-1}{24}$ & $\sfrac{-1}{24}$ \\
   $C_5$ & $\sfrac{107}{6}$ & $\sfrac{73}{3}$ & $\sfrac{191}{6}$ & $\sfrac{121}{3}$ & $\sfrac{299}{6}$ \\
   $C_6$ & -6 & -7 & -8 & -9 & -10 \\
   $C_7$ & $\sfrac{1}{2}$ & $\sfrac{1}{2}$ & $\sfrac{1}{2}$ & $\sfrac{1}{2}$ & $\sfrac{1}{2}$ \\
   $C_8$ & -3 & $\sfrac{-7}{2}$ & -4 & $\sfrac{-9}{2}$ & -5 \\
   $C_9$ & $\sfrac{1}{2}$ & $\sfrac{1}{2}$ & $\sfrac{1}{2}$ & $\sfrac{1}{2}$ & $\sfrac{1}{2}$ \\
   $C_{10}$ & $\sfrac{1}{6}$ & $\sfrac{1}{6}$ & $\sfrac{1}{6}$ & $\sfrac{1}{6}$ & $\sfrac{1}{6}$ \\
   $C_{11}$ & $\sfrac{11}{2}$ & $\sfrac{13}{2}$ & $\sfrac{15}{2}$ & $\sfrac{17}{2}$ & $\sfrac{19}{2}$ \\
   $C_{12}$ & -1 & -1 & -1 & -1 & -1 \\
   $C_{13}$ & -1 & -1 & -1 & -1 & -1 \\
   $C_{14}$ & $\sfrac{-1}{2}$ & $\sfrac{-1}{2}$ & $\sfrac{-1}{2}$ & $\sfrac{-1}{2}$ & $\sfrac{-1}{2}$ \\
   $C_{15}$ & 1 & 1 & 1 & 1 & 1 \\
   \hline
\end{longtable}

Solving for several of the coefficients is trivial, as they remain constant regardless of the order $p$. The remaining coefficients change with $p$; however they can be obtained through a simple curve-fitting procedure. Following this approach, the final nodal index formulation in four dimensions is given by
\begin{align} \label{4D_Index}
    m &=1+\left(\frac{2p^3+15p^2+35p+25}{12}\right)i+\left(\frac{-6p^2-30p-35}{24}\right)i^2 \notag \\[1.0ex]
    &+\left(\frac{2p+5}{12}\right)i^3-\left(\frac{1}{24}\right)i^4+\left(\frac{3p^2+12p+11}{6}\right)j+(-p-2)ij+\left(\frac{1}{2}\right)i^2j \notag \\[1.0ex]
    &+\left(\frac{-p-2}{2}\right)j^2+\left(\frac{1}{2}\right)ij^2+\left(\frac{1}{6}\right)j^3+\left(\frac{2p+3}{2}\right)k-ik-jk-\left(\frac{1}{2}\right)k^2+q, \notag \\[1.0ex]  &(i,j,k,q)\geq 0, \quad i+j+k+q\leq p. 
\end{align}
In order to minimize computational roundoff error, the authors find it useful to premultiply each term in Eq.~\eqref{4D_Index} by a factor of 24 and then multiply the result of the entire equation by $\sfrac{1}{24}$ after summing terms.

\subsection{Lebesgue Constant} \label{lebesgue_sec}

With the indexing function derived, all the preliminary pieces are now in place to generate the nodal basis functions on the pentatope. It remains for us to establish a metric for assessing the validity of the interpolation point sets. This assessment is performed by calculating the Lebesgue constant $\Lambda$ for each point set $\left\{\bm{r}_i \right\}$
\begin{align}
    \Lambda = \max_{\bm{r}} \sum_{i=1}^{N_p} \left| l_i(\bm{r})\right|.
\end{align}
Minimization of the Lebesgue constant is proven to minimize the $L^{\infty}$ error as follows
\begin{align}
    ||u-u_h||_{\infty} = ||u-u^*+u^*-u_h||_{\infty} \leq ||u-u^*||_{\infty} + ||u^*-u_h||_{\infty} \leq (1+\Lambda)||u-u^*||_{\infty}, 
\end{align}
where $u^*$ is the best polynomial approximation of the exact solution $u$.

\section{Methodology} \label{method_sec}

Optimization of the Lebesgue constant for a given point set is performed using warping and blending functions to adjust interpolation point locations. These functions are generated sequentially, with the four-dimensional functions dependent on the three-dimensional ones—which are themselves dependent on the two-dimensional functions and so on. Here, the process of generating these warping and blending functions is detailed, beginning with the one-dimensional case and extending up to four dimensions.

\subsection{The Line Segment}\label{Segment}

Consider the reference line segment, which is defined as the shortest distance between the following two vertices
\begin{align}
    \bm{v}_{1} &= -1, \qquad \bm{v}_{2} = 1.
\end{align}
By construction, the reference line segment is of length 2, and the sole edge is given by the subset consisting of the vertices
\begin{align*}
    E_{1} &= \left\{\bm{v}_{1}, \bm{v}_{2} \right\}.
\end{align*}
The simplest discretization of the reference line is that of equally-spaced nodes
\begin{align*}
    r^e_i = -1+\frac{2i}{p}, \qquad i \in [0,...,p],
\end{align*}
However, this discretization can result in the onset of the infamous Runge phenomenon for large values of $p$. As such, it is often advantageous to adopt a `warped' discretization using Legendre-Gauss-Lobatto point locations, which are optimized to maximize the determinant of the associated Vandermonde matrix. These points will be referred to as $r_i^{LGL}$. These two nodal discretizations can be connected via
\begin{align*}
    w(r) = \sum_{i=1}^{p+1}(r_i^{LGL}-r_i^e)\ell_i^e(r),
\end{align*}
where $\ell_i^e$ are the Lagrange basis functions derived for $r_i^e$.

The warping function (above) can be modified in order to obtain
\begin{align}
    \widetilde{w}(r) = \frac{w(r)}{1-r^2}. \label{warp_mod}
\end{align}
This function will be used frequently during the subsequent sections.

\subsection{The Triangle}\label{Triangle}

Consider the reference triangle, which is defined as the convex hull of the following three vertices
\begin{align}
    \bm{v}_{1} = \left(-1, -\frac{1}{\sqrt{3}} \right), \qquad \bm{v}_{2} = \left(1,-\frac{1}{\sqrt{3}} \right), \qquad \bm{v}_{3} = \left(0, \frac{2}{\sqrt{3}} \right).
\end{align}
All points within the triangle interior can be represented as follows
\begin{align*}
    \bm{r} = \lambda^2 \bm{v}_1 + \lambda^3 \bm{v}_2 + \lambda^1 \bm{v}_3, 
\end{align*}
where
\begin{align*}
    0\leq \lambda^1, \lambda^2, \lambda^3 \leq 1, \qquad \lambda^1 + \lambda^2 + \lambda^3 = 1.
\end{align*}
By construction, the triangle is equilateral and all of the edges are of length 2. The edges are given by
\begin{align*}
    E_{1} &= \left\{\bm{v}_{1}, \bm{v}_{2} \right\}, \qquad
    E_{2} = \left\{\bm{v}_{2}, \bm{v}_{3} \right\}, \qquad 
    E_{3} = \left\{\bm{v}_{1}, \bm{v}_{3} \right\}.
\end{align*}
A warping function for each edge is derived from Eq.~\eqref{warp_mod} as follows
\begin{align*}
    \bm{w}^1(\lambda^1, \lambda^2, \lambda^3) &= \widetilde{w}(\lambda^3-\lambda^2) \bm{t}^{1}, \\[1.0ex]
    \bm{w}^2(\lambda^1, \lambda^2, \lambda^3) &= \widetilde{w}(\lambda^1-\lambda^3) \bm{t}^{2}, \\[1.0ex]
    \bm{w}^3(\lambda^1, \lambda^2, \lambda^3) &= \widetilde{w}(\lambda^2-\lambda^1) \bm{t}^{3}.
\end{align*}
Here, the $\bm{t}^{f}$'s are the vectors associated with each edge
\begin{align*}
    \bm{t}^{1} &= \left(1,0\right), \\[1.0ex]
    \bm{t}^{2} &= \left(-\frac{1}{2}, \frac{\sqrt{3}}{2} \right), \\[1.0ex]
    \bm{t}^{3} &= \left(-\frac{1}{2}, -\frac{\sqrt{3}}{2} \right).
\end{align*}
In addition, we can define the following blending functions for each edge
\begin{align*}
    b^1(\lambda^1,\lambda^2,\lambda^3) = 4 \lambda^3 \lambda^2, \\[1.0ex]
    b^2(\lambda^1,\lambda^2,\lambda^3) = 4 \lambda^3 \lambda^1, \\[1.0ex]
    b^3(\lambda^1,\lambda^2,\lambda^3) = 4 \lambda^2 \lambda^1.
\end{align*}
Next, we can define the following warping/blending function on the entire triangle
\begin{align}
     \bm{g}\left(\lambda^1, \lambda^2, \lambda^3\right) &= \left(1 + \left(\alpha \lambda^1\right)^2\right) b^1 \bm{w}^1 + \left(1 + \left(\alpha \lambda^2\right)^2\right) b^2 \bm{w}^2 + \left(1 + \left(\alpha \lambda^3\right)^2\right) b^3 \bm{w}^3,
\end{align}
where $\alpha$ is an optimization parameter.
Generally speaking, our objective is to shift the locations of equidistant points on the triangle by using the warping/blending function $\bm{g}\left(\lambda^1, \lambda^2, \lambda^3\right)$, and thereby obtain good interpolation points. The locations of the equidistant points are given by
\begin{align*}
    (i,j)\geq 0, \quad i+j\leq p: \qquad \left(\lambda^{1}_{e}, \lambda^{3}_{e} \right) = \left(\frac{i}{p}, \frac{j}{p} \right), \qquad \lambda^{2}_{e} = 1-\lambda^{1}_{e} -\lambda^{3}_{e},
\end{align*}
in barycentric coordinates, or
\begin{align*}
    \bm{r}_{e} = \lambda^{2}_{e} \bm{v}_1 + \lambda^{3}_{e} \bm{v}_2 + \lambda^{1}_{e} \bm{v}_3,
\end{align*}
in reference coordinates. We obtain the final shifted point locations by simply adding the warping/blending function to the equidistant point locations as follows
\begin{align}
    \bm{r}_{s} = \bm{r}_{e} + \bm{g}\left(\lambda^{1}_{e}, \lambda^{2}_{e}, \lambda^{3}_{e} \right).
\end{align}
If the warping/blending function has been formulated correctly, then the shifted points $\bm{r}_{s}$ will act as good-quality interpolation points.

\subsection{The Tetrahedron}\label{Tetrahedron}

\begin{figure}[h!]
\centering
    \includegraphics[width = 0.6\linewidth]{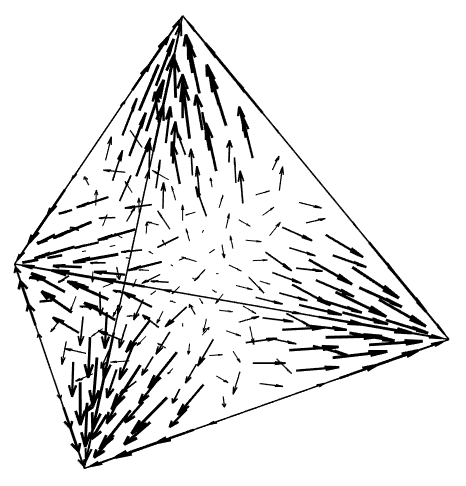}
    \caption{Effect of warping and blending functions on points on the  reference tetrahedron. Note that points located at vertices, midpoints, and barycenters undergo zero warping and/or blending.}
    \label{tet_warp}
\end{figure}

Consider the reference tetrahedron, which is defined as the convex hull of the following four vertices
\begin{align}
    \nonumber \bm{v}_{1} &= \left(-1, -\frac{1}{\sqrt{3}}, -\frac{1}{\sqrt{6}} \right), \qquad \bm{v}_{2} = \left(1,-\frac{1}{\sqrt{3}}, -\frac{1}{\sqrt{6}} \right), \\[1.0ex]
    \bm{v}_{3} &= \left(0, \frac{2}{\sqrt{3}}, - \frac{1}{\sqrt{6}} \right), \qquad \qquad \bm{v}_{4} = \left(0,0,\frac{3}{\sqrt{6}} \right).
\end{align}
All points within the tetrahedron interior can be represented as follows
\begin{align*}
    \bm{r} = \lambda^3 \bm{v}_1 + \lambda^4 \bm{v}_2 + \lambda^2 \bm{v}_3 + \lambda^1 \bm{v}_4, 
\end{align*}
where
\begin{align*}
    0\leq \lambda^1, \lambda^2, \lambda^3, \lambda^4 \leq 1, \qquad \lambda^1 + \lambda^2 + \lambda^3 + \lambda^4 = 1.
\end{align*}
The tetrahedron has a total of four equilateral triangular faces (each with edge length 2), given by the following subsets of vertices
\begin{align*}
    F_{1} &= \left\{\bm{v}_{1}, \bm{v}_{2}, \bm{v}_{3} \right\}, \qquad
    F_{2} = \left\{\bm{v}_{1}, \bm{v}_{2}, \bm{v}_{4} \right\}, \\[1.0ex] 
    F_{3} &= \left\{ \bm{v}_{2}, \bm{v}_{3}, \bm{v}_{4} \right\}, \qquad
    F_{4} = \left\{\bm{v}_{1}, \bm{v}_{3}, \bm{v}_{4} \right\}.
\end{align*}
Let us now define the warping functions for each face as follows
\begin{align*}
     \bm{w}^{1} &= \bm{g}_{1} \left(\lambda^2, \lambda^3, \lambda^4 \right) \bm{t}^{1,1} + \bm{g}_{2} \left(\lambda^2, \lambda^3, \lambda^4 \right) \bm{t}^{1,2},  \\[1.0ex]
    \bm{w}^{2} &= \bm{g}_{1} \left(\lambda^1, \lambda^3, \lambda^4 \right) \bm{t}^{2,1} + \bm{g}_{2} \left(\lambda^1, \lambda^3, \lambda^4 \right) \bm{t}^{2,2}, \\[1.0ex]
    \bm{w}^{3} &= \bm{g}_{1} \left( \lambda^1, \lambda^4, \lambda^2\right) \bm{t}^{3,1} + \bm{g}_{2} \left(\lambda^1, \lambda^4, \lambda^2 \right) \bm{t}^{3,2}, \\[1.0ex]
    \bm{w}^{4} &= \bm{g}_{1} \left(\lambda^1, \lambda^3, \lambda^2 \right) \bm{t}^{4,1} + \bm{g}_{2} \left( \lambda^1, \lambda^3, \lambda^2 \right) \bm{t}^{4,2}.
\end{align*}
Here, $\bm{t}^{f,1}$ and $\bm{t}^{f,2}$ are two orthonormal vectors associated with each face
\begin{align*}
    \bm{t}^{1,1} &= \left(1,0,0 \right), \quad     \bm{t}^{1,2} = \left(0,1,0 \right), \\[1.0ex]
    \bm{t}^{2,1} &= \left(1,0,0 \right), \quad     \bm{t}^{2,2} = \left(0, \frac{1}{3}, \frac{4}{3\sqrt{2}} \right), \\[1.0ex]
    \bm{t}^{3,1} &= \left(-\frac{1}{2},\frac{3}{2\sqrt{3}},0 \right), \quad     \bm{t}^{3,2} = \left(-\frac{1}{2\sqrt{3}},-\frac{1}{6},\frac{4}{3\sqrt{2}} \right), \\[1.0ex]
    \bm{t}^{4,1} &=   \left(\frac{1}{2},\frac{3}{2\sqrt{3}},0 \right), \quad     \bm{t}^{4,2} = \left(\frac{1}{2\sqrt{3}},-\frac{1}{6}, \frac{4}{3\sqrt{2}}\right).
\end{align*}
Next, we can define the following blending functions for each face
\begin{align*}
    b^{1}(\lambda^1,\lambda^2,\lambda^3,\lambda^4) &= \left( \frac{2\lambda^2}{2\lambda^2 + \lambda^1} \right)\left( \frac{2\lambda^3}{2\lambda^3 + \lambda^1} \right)\left( \frac{2\lambda^4}{2\lambda^4 + \lambda^1} \right), \\[1.0ex]
    b^{2}(\lambda^1,\lambda^2,\lambda^3,\lambda^4) &= \left( \frac{2\lambda^1}{2\lambda^1 + \lambda^2} \right)\left( \frac{2\lambda^3}{2\lambda^3 + \lambda^2} \right)\left( \frac{2\lambda^4}{2\lambda^4 + \lambda^2} \right), \\[1.0ex]
    b^{3}(\lambda^1,\lambda^2,\lambda^3,\lambda^4) &= \left( \frac{2\lambda^1}{2\lambda^1 + \lambda^3} \right)\left( \frac{2\lambda^2}{2\lambda^2 + \lambda^3} \right)\left( \frac{2\lambda^4}{2\lambda^4 + \lambda^3} \right), \\[1.0ex]
    b^{4}(\lambda^1,\lambda^2,\lambda^3,\lambda^4) &= \left( \frac{2\lambda^1}{2\lambda^1 + \lambda^4} \right)\left( \frac{2\lambda^2}{2\lambda^2 + \lambda^4} \right)\left( \frac{2\lambda^3}{2\lambda^3 + \lambda^4} \right).
\end{align*}
Furthermore, we can define the following warping/blending function on the entire tetrahedron
\begin{align}
    \nonumber \bm{g}\left(\lambda^1, \lambda^2, \lambda^3, \lambda^4\right) &= \left(1 + \left(\beta \lambda^1\right)^2\right) b^1 \bm{w}^1 + \left(1 + \left(\beta \lambda^2\right)^2\right) b^2 \bm{w}^2  \\[1.0ex]
    &+ \left(1 + \left(\beta \lambda^3\right)^2\right) b^3 \bm{w}^3 + \left(1 + \left(\beta \lambda^4\right)^2\right) b^4 \bm{w}^4,
\end{align}
where $\beta$ is an optimization parameter. Figure~\ref{tet_warp} illustrates the combined impact of the warping and blending functions on the tetrahedron.

As before, we can use the function $\bm{g}\left(\lambda^1, \lambda^2, \lambda^3, \lambda^4 \right)$ to shift the locations of equidistant points on the tetrahedron. The locations of the equidistant points are given by
\begin{align*}
    &(i,j,k)\geq 0, \quad i+j+k\leq p: \qquad \left(\lambda^{1}_{e}, \lambda^{2}_{e}, \lambda^{4}_{e} \right) = \left(\frac{i}{p}, \frac{j}{p}, \frac{k}{p} \right), \\[1.0ex] 
    &\lambda^{3}_{e} = 1-\lambda^{1}_{e} -\lambda^{2}_{e} - \lambda^{4}_{e},
\end{align*}
in barycentric coordinates, or
\begin{align*}
    \bm{r}_{e} = \lambda^{3}_{e} \bm{v}_1 + \lambda^{4}_{e} \bm{v}_2 + \lambda^{2}_{e} \bm{v}_3 + \lambda^{1}_{e} \bm{v}_4,
\end{align*}
in reference coordinates. We obtain the final shifted point locations by adding the warping/blending function to the equidistant point locations as follows
\begin{align}
    \bm{r}_{s} = \bm{r}_{e} + \bm{g}\left(\lambda^{1}_{e}, \lambda^{2}_{e}, \lambda^{3}_{e}, \lambda^{4}_{e} \right).
\end{align}

\subsection{The Pentatope}\label{Pentatope}

Consider the reference pentatope, which is defined as the convex hull of the following five vertices
\begin{align}
    \nonumber \bm{v}_{1} &= \left(-1, -\frac{1}{\sqrt{3}}, -\frac{1}{\sqrt{6}}, - \frac{1}{\sqrt{10}} \right), \qquad \bm{v}_{2} = \left(1,-\frac{1}{\sqrt{3}}, -\frac{1}{\sqrt{6}}, - \frac{1}{\sqrt{10}} \right), \\[1.0ex]
    \nonumber \bm{v}_{3} &= \left(0, \frac{2}{\sqrt{3}}, - \frac{1}{\sqrt{6}}, - \frac{1}{\sqrt{10}} \right), \, \; \quad \qquad \bm{v}_{4} = \left(0,0,\frac{3}{\sqrt{6}}, -\frac{1}{\sqrt{10}} \right), \\[1.0ex]
    \bm{v}_{5} &= \left(0,0,0, \frac{4}{\sqrt{10}} \right).
\end{align}
All points within the pentatope interior can be represented as follows
\begin{align*}
    \bm{r} = \lambda^4 \bm{v}_1 + \lambda^5 \bm{v}_2 + \lambda^3 \bm{v}_3 + \lambda^2 \bm{v}_4 + \lambda^1 \bm{v}_5, 
\end{align*}
where
\begin{align*}
    0\leq \lambda^1, \lambda^2, \lambda^3, \lambda^4, \lambda^5 \leq 1, \qquad \lambda^1 + \lambda^2 + \lambda^3 + \lambda^4 + \lambda^5 = 1.
\end{align*}
By construction, the pentatope is equi-facetal, and all of the edges are of length~2. In addition, the pentatope has a total of five tetrahedral facets, given by the following subsets of vertices
\begin{align*}
    F_{1} &= \left\{\bm{v}_{1}, \bm{v}_{2}, \bm{v}_{3}, \bm{v}_{4} \right\}, \qquad
    F_{2} = \left\{\bm{v}_{1}, \bm{v}_{2}, \bm{v}_{3}, \bm{v}_{5} \right\}, \\[1.0ex] 
    F_{3} &= \left\{ \bm{v}_{1}, \bm{v}_{2}, \bm{v}_{4}, \bm{v}_{5} \right\}, \qquad
    F_{4} = \left\{\bm{v}_{2}, \bm{v}_{3}, \bm{v}_{4}, \bm{v}_{5} \right\}, \\[1.0ex] 
    F_{5} &= \left\{\bm{v}_{3}, \bm{v}_{1}, \bm{v}_{4}, \bm{v}_{5}, \right\}.
\end{align*}
We can now define vector warping functions for each facet as follows
\begin{align*}
    \bm{w}^{1} &= \bm{g}_{1} \left(\lambda^2, \lambda^3, \lambda^4, \lambda^5 \right) \bm{t}^{1,1} + \bm{g}_{2} \left(\lambda^2, \lambda^3, \lambda^4, \lambda^5 \right) \bm{t}^{1,2} + \bm{g}_{3} \left(\lambda^2, \lambda^3, \lambda^4, \lambda^5 \right) \bm{t}^{1,3},  \\[1.0ex]
    \bm{w}^{2} &= \bm{g}_{1} \left( \lambda^1,\lambda^3,\lambda^4,\lambda^5 \right) \bm{t}^{2,1} + \bm{g}_{2} \left(\lambda^1,\lambda^3,\lambda^4,\lambda^5 \right) \bm{t}^{2,2} + \bm{g}_{3} \left(\lambda^1,\lambda^3,\lambda^4,\lambda^5 \right) \bm{t}^{2,3}, \\[1.0ex]
    \bm{w}^{3} &= \bm{g}_{1} \left(\lambda^1,\lambda^2,\lambda^4,\lambda^5 \right) \bm{t}^{3,1} + \bm{g}_{2} \left(\lambda^1,\lambda^2,\lambda^4,\lambda^5 \right) \bm{t}^{3,2} + \bm{g}_{3} \left(\lambda^1,\lambda^2,\lambda^4,\lambda^5 \right) \bm{t}^{3,3}, \\[1.0ex]
    \bm{w}^{4} &= \bm{g}_{1} \left( \lambda^1,\lambda^2,\lambda^5,\lambda^3 \right) \bm{t}^{4,1} + \bm{g}_{2} \left( \lambda^1,\lambda^2,\lambda^5,\lambda^3 \right) \bm{t}^{4,2} + \bm{g}_{3} \left( \lambda^1,\lambda^2,\lambda^5,\lambda^3 \right) \bm{t}^{4,3}, \\[1.0ex]
    \bm{w}^{5} &= \bm{g}_{1} \left(\lambda^1,\lambda^2,\lambda^4,\lambda^3 \right) \bm{t}^{5,1} + \bm{g}_{2} \left(\lambda^1,\lambda^2,\lambda^4,\lambda^3 \right) \bm{t}^{5,2} + \bm{g}_{3} \left(\lambda^1,\lambda^2,\lambda^4,\lambda^3 \right) \bm{t}^{5,3}.
\end{align*}
Here, $\bm{t}^{f,1}$, $\bm{t}^{f,2}$, and $\bm{t}^{f,3}$ are three orthonormal vectors associated with each facet
\begin{align*}
    \bm{t}^{1,1} &= \left(1,0,0,0\right), \quad \bm{t}^{1,2} = \left(0,1,0,0\right), \quad \bm{t}^{1,3} = \left(0,0,1,0\right), \\[1.0ex]
    \bm{t}^{2,1} &= \left(1,0,0,0\right), \quad \bm{t}^{2,2} = \left(0,1,0,0\right), \quad \bm{t}^{2,3} = \left(0,0,\frac{1}{4}, \frac{\sqrt{15}}{4}\right), \\[1.0ex]
    \bm{t}^{3,1} &= \left(1,0,0,0\right), \quad \bm{t}^{3,2} = \left(0,\frac{1}{3},\frac{4}{3\sqrt{2}},0\right), \quad \bm{t}^{3,3} = \left(0,\frac{1}{\sqrt{18}},-\frac{1}{12},\frac{\sqrt{15}}{4} \right), \\[1.0ex]
    \bm{t}^{4,1} &= \left(-\frac{1}{2},\frac{3}{2\sqrt{3}},0,0\right), \quad \bm{t}^{4,2} = \left(-\frac{1}{2\sqrt{3}},-\frac{1}{6},\frac{4}{3\sqrt{2}},0\right), \quad \bm{t}^{4,3} = \left(-\frac{\sqrt{6}}{12},-\frac{\sqrt{2}}{12},-\frac{1}{12},\frac{\sqrt{15}}{4} \right), \\[1.0ex]
    \bm{t}^{5,1} &= \left(\frac{1}{2},\frac{3}{2\sqrt{3}},0,0\right), \quad \bm{t}^{5,2} = \left(\frac{1}{2\sqrt{3}},-\frac{1}{6},\frac{4}{3\sqrt{2}},0\right), \quad \bm{t}^{5,3} = \left(\frac{\sqrt{6}}{12},-\frac{\sqrt{2}}{12},-\frac{1}{12},\frac{\sqrt{15}}{4} \right).
\end{align*}
In addition, we can define the following blending functions for each facet
\begin{align*}
    b^{1}(\lambda^1,\lambda^2,\lambda^3,\lambda^4,\lambda^5) &= \left( \frac{2\lambda^2}{2\lambda^2 + \lambda^1} \right)\left( \frac{2\lambda^3}{2\lambda^3 + \lambda^1} \right)\left( \frac{2\lambda^4}{2\lambda^4 + \lambda^1} \right)\left( \frac{2\lambda^5}{2\lambda^5 + \lambda^1} \right), \\[1.0ex]
    b^{2}(\lambda^1,\lambda^2,\lambda^3,\lambda^4,\lambda^5) &= \left( \frac{2\lambda^1}{2\lambda^1 + \lambda^2} \right)\left( \frac{2\lambda^3}{2\lambda^3 + \lambda^2} \right)\left( \frac{2\lambda^4}{2\lambda^4 + \lambda^2} \right)\left( \frac{2\lambda^5}{2\lambda^5 + \lambda^2} \right), \\[1.0ex]
    b^{3}(\lambda^1,\lambda^2,\lambda^3,\lambda^4,\lambda^5) &= \left( \frac{2\lambda^1}{2\lambda^1 + \lambda^3} \right)\left( \frac{2\lambda^2}{2\lambda^2 + \lambda^3} \right)\left( \frac{2\lambda^4}{2\lambda^4 + \lambda^3} \right)\left( \frac{2\lambda^5}{2\lambda^5 + \lambda^3} \right), \\[1.0ex]
    b^{4}(\lambda^1,\lambda^2,\lambda^3,\lambda^4,\lambda^5) &= \left( \frac{2\lambda^1}{2\lambda^1 + \lambda^4} \right)\left( \frac{2\lambda^2}{2\lambda^2 + \lambda^4} \right)\left( \frac{2\lambda^3}{2\lambda^3 + \lambda^4} \right)\left( \frac{2\lambda^5}{2\lambda^5 + \lambda^4} \right), \\[1.0ex]
    b^{5}(\lambda^1,\lambda^2,\lambda^3,\lambda^4,\lambda^5) &= \left( \frac{2\lambda^1}{2\lambda^1 + \lambda^5} \right)\left( \frac{2\lambda^2}{2\lambda^2 + \lambda^5} \right)\left( \frac{2\lambda^3}{2\lambda^3 + \lambda^5} \right)\left( \frac{2\lambda^4}{2\lambda^4 + \lambda^5} \right).
\end{align*}
Next, we can define the following warping/blending function on the entire pentatope
\begin{align} \label{penta_warpblend}
    \nonumber \bm{g}\left(\lambda^1, \lambda^2, \lambda^3, \lambda^4, \lambda^5 \right) &= \left(1 + \left(\gamma \lambda^1\right)^2\right) b^1 \bm{w}^1 + \left(1 + \left(\gamma \lambda^2\right)^2\right) b^2 \bm{w}^2 + \left(1 + \left(\gamma \lambda^3\right)^2\right) b^3 \bm{w}^3 \\[1.0ex]
    &+ \left(1 + \left(\gamma \lambda^4\right)^2\right) b^4 \bm{w}^4 + \left(1 + \left(\gamma \lambda^5\right)^2\right) b^5 \bm{w}^5,
\end{align}
where $\gamma$ is an optimization parameter. As before, we can use the function $\bm{g}\left(\lambda^1, \lambda^2, \lambda^3, \lambda^4, \lambda^5 \right)$ to shift the locations of equidistant points on the pentatope. The locations of the equidistant points are given by
\begin{align*}
    (i,j,k,m)\geq 0, \quad i+j+k+m\leq p&: \qquad \left(\lambda^{1}_{e}, \lambda^{2}_{e}, \lambda^{3}_{e}, \lambda^{5}_{e} \right) = \left(\frac{i}{p}, \frac{j}{p}, \frac{k}{p}, \frac{m}{p} \right), \\[1.0ex] 
    \lambda^{4}_{e} = 1-\lambda^{1}_{e} &-\lambda^{2}_{e} - \lambda^{3}_{e} - \lambda^{5}_{e},
\end{align*}
in barycentric coordinates, or
\begin{align*}
   \bm{r}_{e} = \lambda^{4}_{e} \bm{v}_1 + \lambda^{5}_{e} \bm{v}_2 + \lambda^{3}_{e} \bm{v}_3 + \lambda^{2}_{e} \bm{v}_4 + \lambda^{1}_{e} \bm{v}_5,
\end{align*}
in reference coordinates. We obtain the final shifted point locations by adding the warping/blending function to the equidistant point locations as follows
\begin{align}
    \bm{r}_{s} = \bm{r}_{e} + \bm{g}\left(\lambda^{1}_{e}, \lambda^{2}_{e}, \lambda^{3}_{e}, \lambda^{4}_{e}, \lambda^{5}_{e} \right).
\end{align}

\renewcommand*{\thefootnote}{\fnsymbol{footnote}}

\section{Results} \label{results_sec}

In the interest of expediting computation we assumed that our three optimization parameters $\alpha, \beta,$ and $\gamma$ were all equivalent, thereby leading to a single-parameter optimization scheme with $\alpha = \beta = \gamma$. The resultant single optimization parameter is henceforth referred to simply as $\alpha$. We note that when $\alpha = 0$, we counterintuitively still obtain nontrivial warping and blending of our equidistant point sets, as evidenced by Eq.~\eqref{penta_warpblend}.

In order to find the optimal value of $\alpha$ for a given polynomial order $p$, a simple bisection algorithm was created. The initial range for this algorithm was $\alpha = [0, 3]$, as preliminary computational results showed that for $\alpha$ values greater than 3, the Lebesgue constant of the shifted nodes tended to be higher than that of the equispaced nodes. The resultant optimized values of $\alpha$ for polynomial orders $p = 1$ to $p = 10$ are tabulated in Table~\ref{alpha_table_summary}, alongside the associated Lebesgue constants $\Lambda$. It is clear from the table that in virtually all cases the optimized points outperform the equispaced points, i.e.~they have smaller Lebesgue constants. The only exception is the $p = 3$ case, for which we obtain a slightly worse Lebesgue constant relative to the equispaced points. 

The Lebesgue constants in Table~\ref{alpha_table_summary} were obtained using a naive sampling approach. Absolute values of each nodal basis function were calculated on a uniform grid within the pentatope (this uniform grid was independent of the equispaced points used for interpolation purposes). Due to rising computational cost, higher polynomial order cases were carried out using a progressively coarser uniform grid. The grid spacing used for each polynomial order is given in the penultimate column of Table~\ref{alpha_table_summary}. Conveniently, previous researchers have computed Lebesgue constants for equispaced interpolation points up to dimension 6 and $p = 10$~\cite{jimenezrefining}. As expected, our equispaced results are in close agreement with theirs for~$d = 4$.  

\begin{figure}[h!]
\centering
\includegraphics[width=10cm,trim={1cm 1cm 1cm 0.5cm},clip]{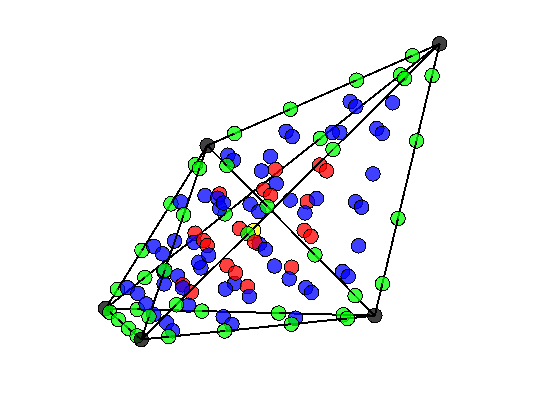}
\caption{Optimized interpolation points on a projected version of the pentatope for $p = 5$. Black points appear at the vertices, green points appear on the edges, blue points appear on the triangular faces, red points appear on the tetrahedral facets, and yellow points appear in the interior of the pentatope.}
\label{p5_fig_full}
\end{figure}

\begin{figure}[h!]
\centering
\includegraphics[width=10cm,trim={1cm 1cm 1cm 0.5cm},clip]{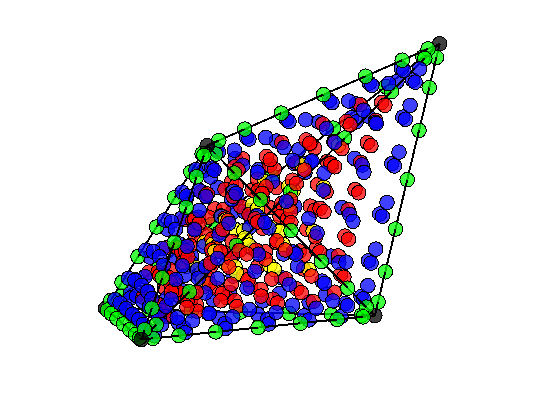}
\caption{Optimized interpolation points on a projected version of the pentatope for $p = 8$. Black points appear at the vertices, green points appear on the edges, blue points appear on the triangular faces, red points appear on the tetrahedral facets, and yellow points appear in the interior of the pentatope.}
\label{p8_fig_full}
\end{figure}

\begin{figure}[h!]
\includegraphics[width=7cm,trim={1cm 1cm 1cm 0},clip]{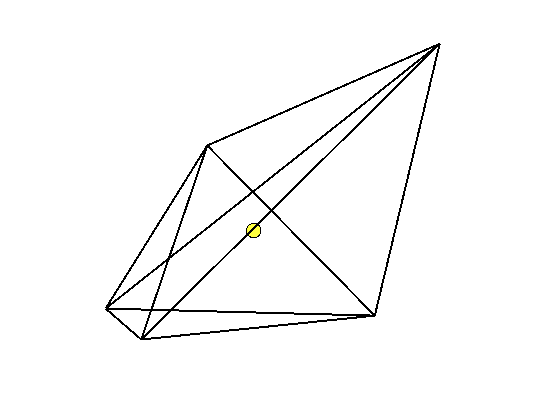}
\includegraphics[width=7cm,trim={1cm 1cm 1cm 0},clip]{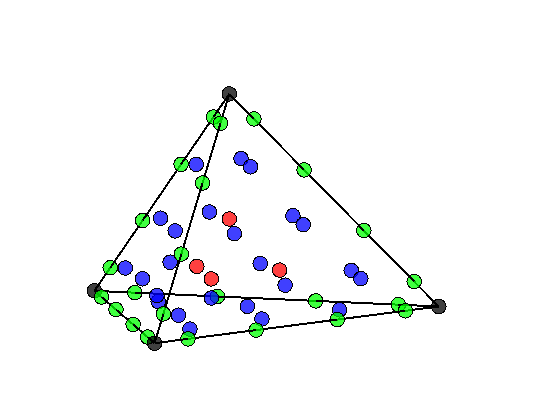}
\includegraphics[width=7cm,trim={1cm 1cm 1cm 0},clip]{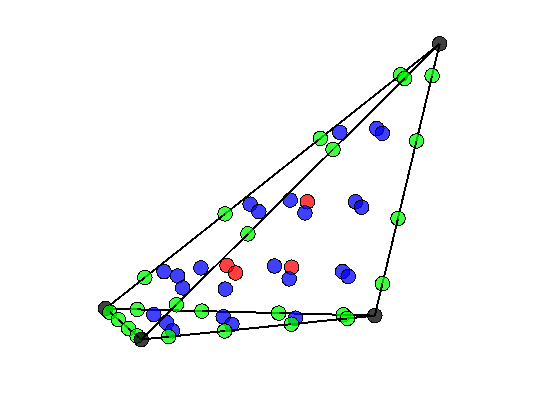}
\includegraphics[width=7cm,trim={1cm 1cm 1cm 0},clip]{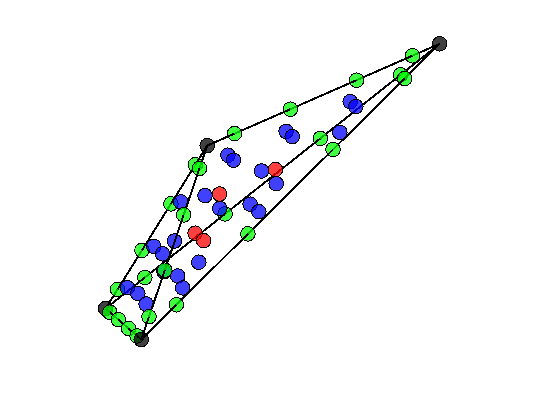}
\includegraphics[width=7cm,trim={1cm 1cm 1cm 0},clip]{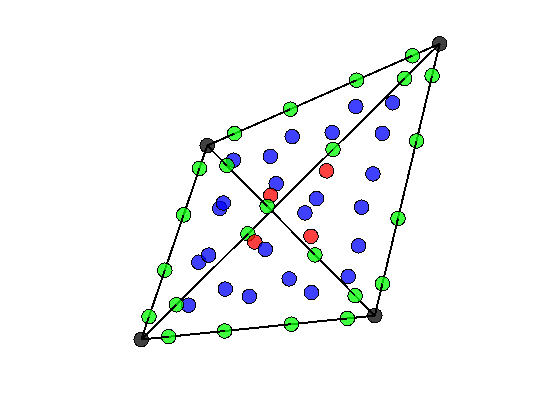}
\includegraphics[width=7cm,trim={1cm 1cm 1cm 0},clip]{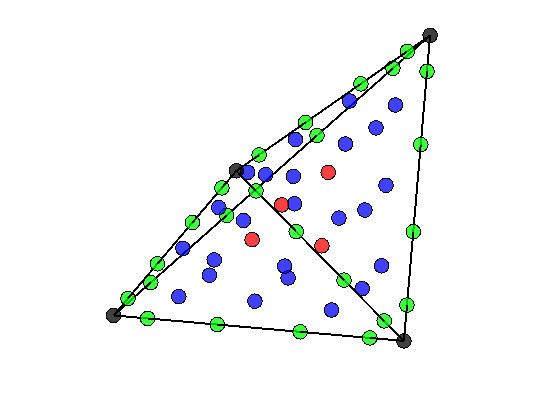}
\caption{Optimized interpolation points on a projected version of the pentatope for $p = 5$, expanded view. Pentatope interior (top left), facet 1 (top right), facet 2 (middle left), facet 3 (middle right), facet 4 (bottom left), and facet 5 (bottom right).  Black points appear at the vertices, green points appear on the edges, blue points appear on the triangular faces, red points appear on the tetrahedral facets, and yellow points appear in the interior of the pentatope.}
\label{p5_fig_facet}
\end{figure}

\begin{figure}[h!]
\includegraphics[width=7cm,trim={1cm 1cm 1cm 0},clip]{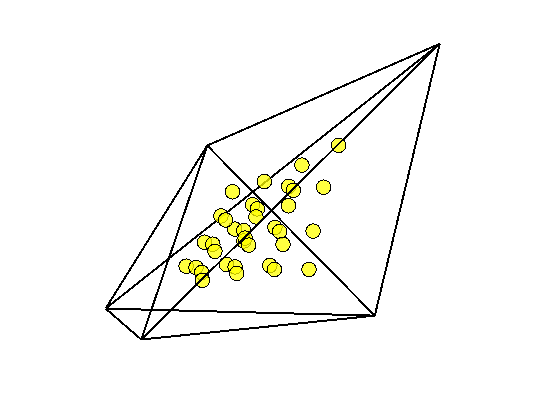}
\includegraphics[width=7cm,trim={1cm 1cm 1cm 0},clip]{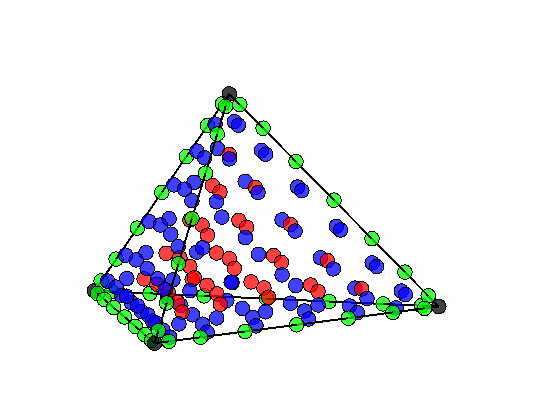}
\includegraphics[width=7cm,trim={1cm 1cm 1cm 0},clip]{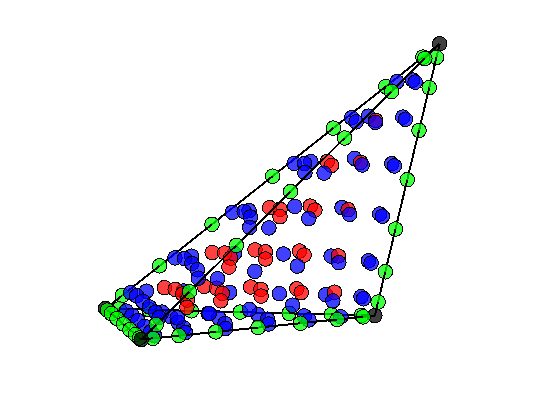}
\includegraphics[width=7cm,trim={1cm 1cm 1cm 0},clip]{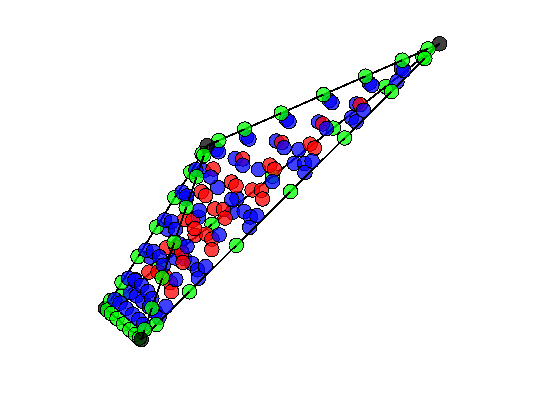}
\includegraphics[width=7cm,trim={1cm 1cm 1cm 0},clip]{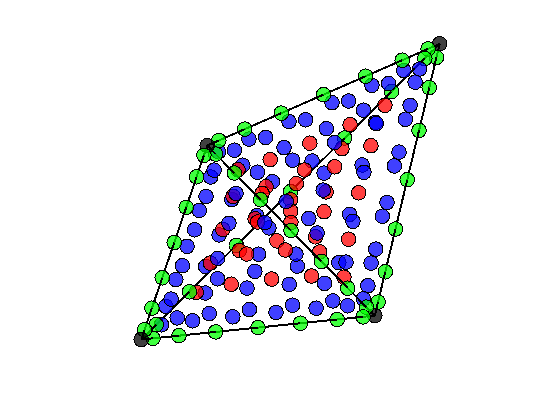}
\includegraphics[width=7cm,trim={1cm 1cm 1cm 0},clip]{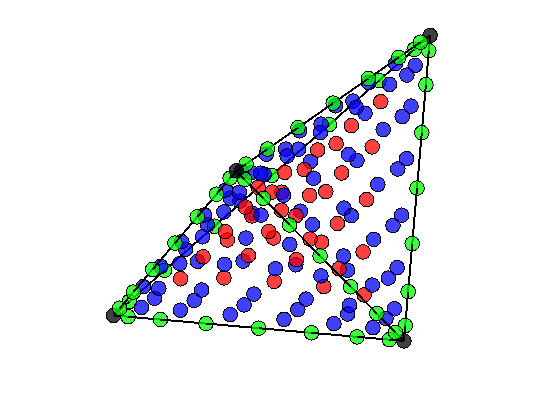}
\caption{Optimized interpolation points on a projected version of the pentatope for $p = 8$, expanded view. Pentatope interior (top left), facet 1 (top right), facet 2 (middle left), facet 3 (middle right), facet 4 (bottom left), and facet 5 (bottom right).  Black points appear at the vertices, green points appear on the edges, blue points appear on the triangular faces, red points appear on the tetrahedral facets, and yellow points appear in the interior of the pentatope.}
\label{p8_fig_facet}
\end{figure}

In addition, the geometric arrangements of the optimized points from $p=5$ to $p=8$ are illustrated in Figures~\ref{p5_fig_full}—\ref{p8_fig_full}~(compact view) as well as~\ref{p5_fig_facet}—\ref{p8_fig_facet}~(expanded view). Finally, the barycentric coordinates for optimized point sets for $p = 1$ to $p = 6$ appear in~\ref{point_tabs}, and the remaining coordinates for $p = 7$ through $p = 10$ appear in~\cite{Gobel24}.

\begin{table}[h!]
\centering
\caption{Comparison of Lebesgue constants for varying polynomial orders $p$ and node generation parameters $\alpha$ on the reference pentatope}
\begin{tabular}{c c c c c c || c}
    \hline
    \multirow{2}{3ex}{$p$} & \multirow{2}{3ex}{$\alpha$} & $\Lambda$ & $\Lambda$ & $\Lambda$ & Grid & $\Lambda$ \\
    & & Optimized & $\alpha = 0$ & Equidistant & Spacing & Equidistant~\cite{jimenezrefining} \\
    \hline
    1 & 0 & 1.0000 & 1.0000 & 1.0000 & 0.01 & \footref{table_dash} \\
    2 & 0 & 2.2000 & 2.2000 & 2.2000 & 0.01 & \footref{table_dash} \\
    3 & 0 & 4.2000 & 4.2000 & 3.8800 & 0.01 & \footref{table_dash} \\
    4 & 0 & 6.1240 & 6.1240 & 6.2384 & 0.01 & \footref{table_dash} \\
    5 & 0 & 8.6423 & 8.6423 & 10.9171 & 0.01 & \tablefootnote{\label{table_dash}Data not given in reference~\cite{jimenezrefining}.} \\
    6 & 1.5000 & 12.0326 & 12.1297 & 19.1418 & 0.02 & 19.22 \\
    7 & 2.2500 & 17.1032 & 18.0971 & 33.8915 & 0.02 & 34.08 \\
    8 & 1.8890 & 23.9226 & 26.0423 & 60.5048 & 0.02 & 60.86 \\
    9 & 1.5000 & 36.1110 & 39.1431 & 109.4267 & 0.02 & 109.43 \\
    10 & 1.5469 & 53.3404 & 57.7742 & 194.8739 & 0.04 & 198.08 \\
    \hline
\end{tabular}
 \label{alpha_table_summary}
\end{table}

\renewcommand*{\thefootnote}{\arabic{footnote}}

\begin{figure}[h!]
\centering
\includegraphics[width=10cm,trim={1cm 1cm 1cm 0.5cm},clip]{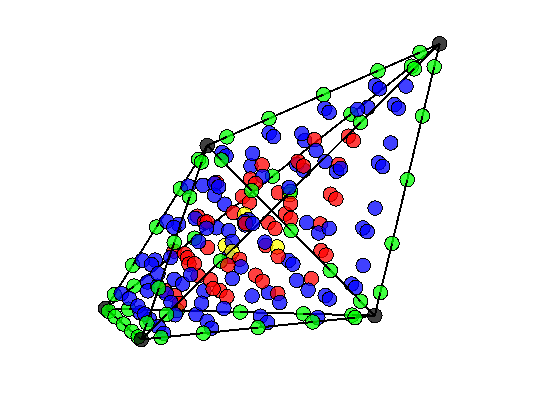}
\caption{Optimized interpolation points on a projected version of the pentatope for $p = 6$. Black points appear at the vertices, green points appear on the edges, blue points appear on the triangular faces, red points appear on the tetrahedral facets, and yellow points appear in the interior of the pentatope.}
\label{p6_fig_full}
\end{figure}

\begin{figure}[h!]
\centering
\includegraphics[width=10cm,trim={1cm 1cm 1cm 0.5cm},clip]{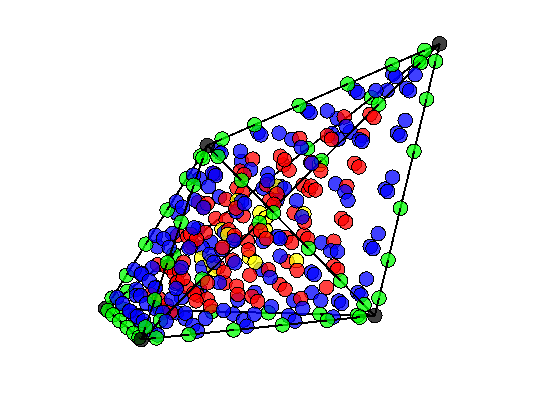}
\caption{Optimized interpolation points on a projected version of the pentatope for $p = 7$. Black points appear at the vertices, green points appear on the edges, blue points appear on the triangular faces, red points appear on the tetrahedral facets, and yellow points appear in the interior of the pentatope.}
\label{p7_fig_full}
\end{figure}

\begin{figure}[h!]
\includegraphics[width=7cm,trim={1cm 1cm 1cm 0},clip]{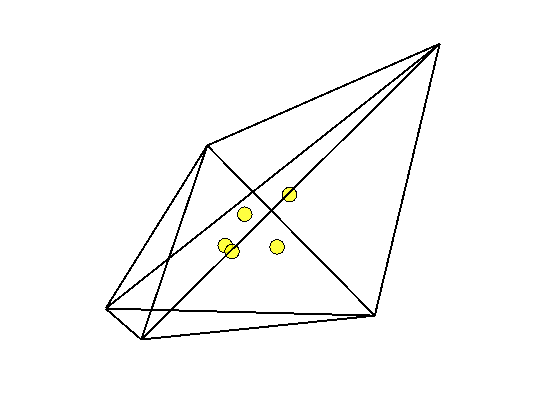}
\includegraphics[width=7cm,trim={1cm 1cm 1cm 0},clip]{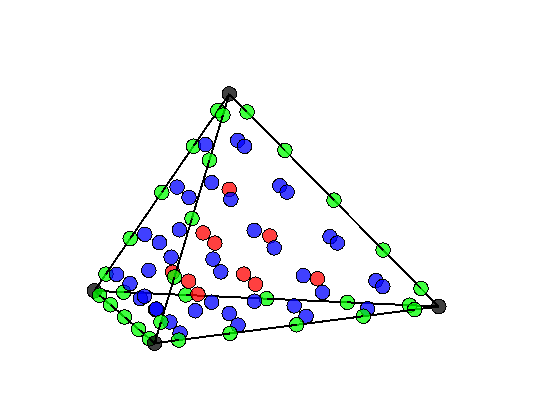}
\includegraphics[width=7cm,trim={1cm 1cm 1cm 0},clip]{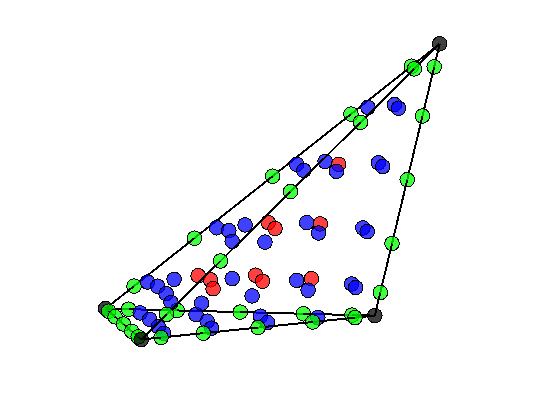}
\includegraphics[width=7cm,trim={1cm 1cm 1cm 0},clip]{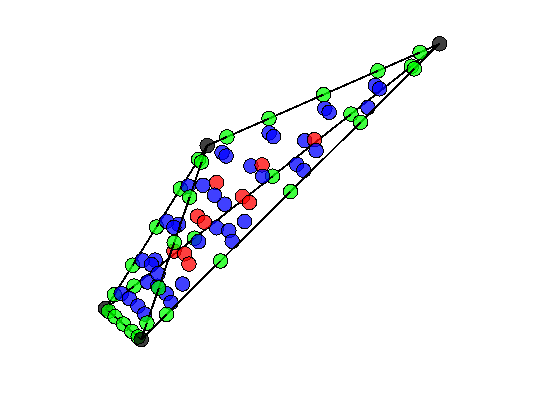}
\includegraphics[width=7cm,trim={1cm 1cm 1cm 0},clip]{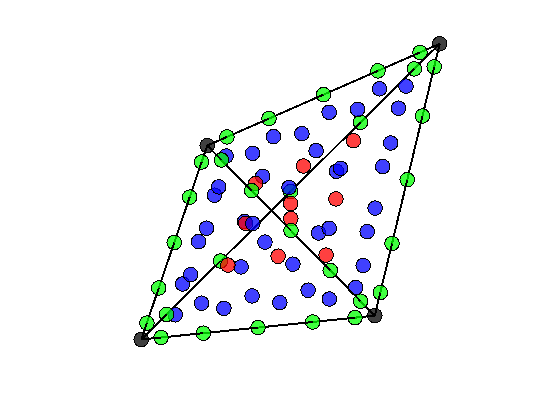}
\includegraphics[width=7cm,trim={1cm 1cm 1cm 0},clip]{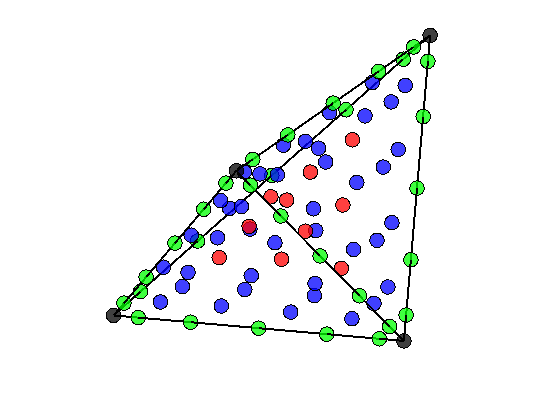}
\caption{Optimized interpolation points on a projected version of the pentatope for $p = 6$, expanded view. Pentatope interior (top left), facet 1 (top right), facet 2 (middle left), facet 3 (middle right), facet 4 (bottom left), and facet 5 (bottom right).  Black points appear at the vertices, green points appear on the edges, blue points appear on the triangular faces, red points appear on the tetrahedral facets, and yellow points appear in the interior of the pentatope.}
\label{p6_fig_facet}
\end{figure}

\begin{figure}[h!]
\includegraphics[width=7cm,trim={1cm 1cm 1cm 0},clip]{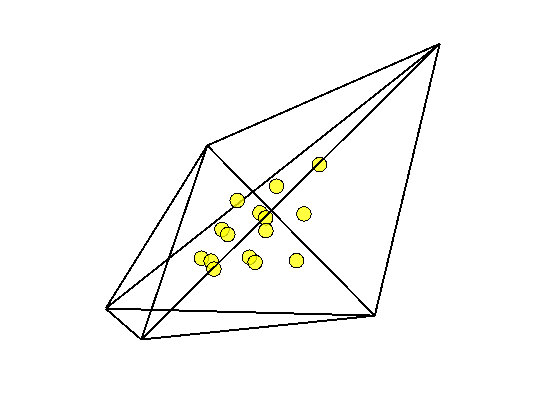}
\includegraphics[width=7cm,trim={1cm 1cm 1cm 0},clip]{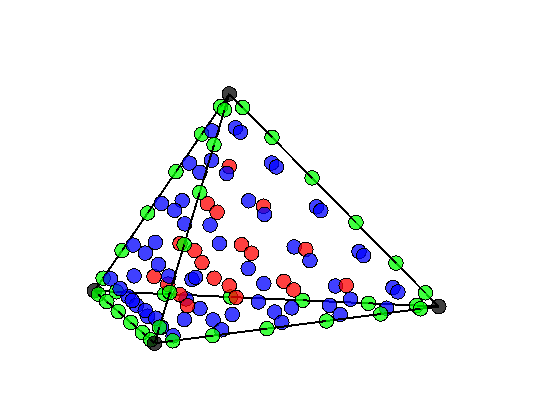}
\includegraphics[width=7cm,trim={1cm 1cm 1cm 0},clip]{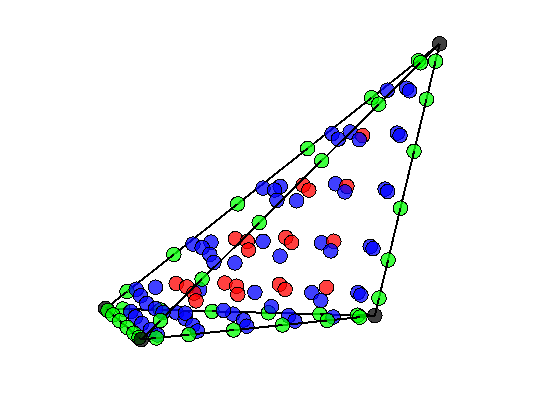}
\includegraphics[width=7cm,trim={1cm 1cm 1cm 0},clip]{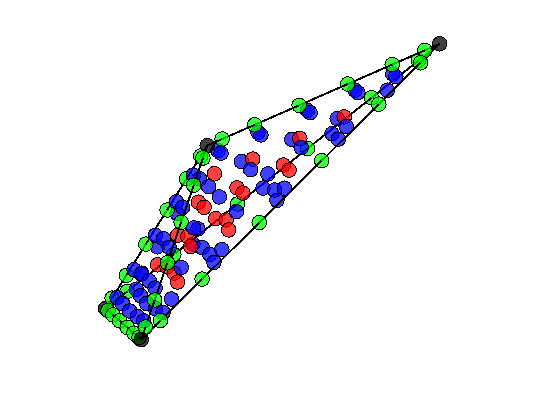}
\includegraphics[width=7cm,trim={1cm 1cm 1cm 0},clip]{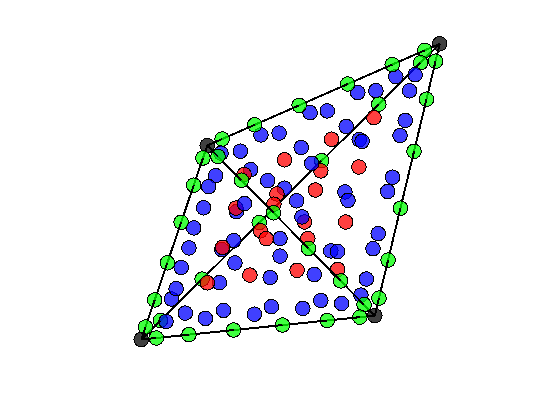}
\includegraphics[width=7cm,trim={1cm 1cm 1cm 0},clip]{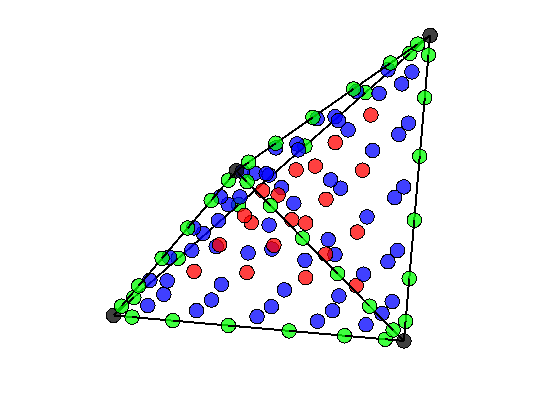}
\caption{Optimized interpolation points on a projected version of the pentatope for $p = 7$, expanded view. Pentatope interior (top left), facet 1 (top right), facet 2 (middle left), facet 3 (middle right), facet 4 (bottom left), and facet 5 (bottom right).  Black points appear at the vertices, green points appear on the edges, blue points appear on the triangular faces, red points appear on the tetrahedral facets, and yellow points appear in the interior of the pentatope.}
\label{p7_fig_facet}
\end{figure}

\pagebreak
\clearpage

\section{Concluding Remarks} \label{conclusion_sec}

In this paper, we generated Lebesgue-optimized interpolation points on four-dimensional simplex elements using an extension of Warburton's approach~\cite{Warburton06} which was originally formulated for problems in two and three dimensions. To this end, we first gave an overview of prior interpolation research on lower dimensional simplices. Then, we constructed the extension of Warburton's warping and blending functions in four dimensions, and discussed how these functions can be applied to the equispaced interpolation point set. Using these functions, we optimized the interpolation point sets around a single parameter $\alpha$ and we found that for polynomial orders $p \leq 5$, the smallest Lebesgue constants occurred for $\alpha = 0$.  This minimal warping and blending—due primarily to the 1D Legendre-Gauss-Lobatto distribution used as a starting point—produces a slightly suboptimal point distribution for $p=3$. However the optimized points began to outperform the equispaced points beginning with $p = 4$. For values of $p > 5$, the optimization parameter $\alpha$ assumed non-zero values. In addition, the gap between the Lebesgue constants of the equispaced points and those of the optimized points rapidly grew. For example, for $p = 10$ the Lebesgue constant of the optimized points is approximately one-fourth of that of the equispaced points.

\subsection{Future Work}

The interpolation point generation procedure detailed in this work can be conceivably extended to higher dimensions, $d > 4$. However, this would require careful construction of the blending/warping functions. In addition, a higher-dimensional extension of the nodal indexing equation would be advantageous. Using the existing indexing equations along with the one derived here, some patterns can be observed regarding coefficients and the terms upon which they act. These patterns will prove useful should a need for the higher-dimensional indexing equations arise.

The simplifications introduced by the present paper can also be modified or omitted entirely. For example, it remains to be shown that equating the optimization parameters of $\alpha, \beta,$ and $\gamma$ leads to the most optimized interpolation point sets. However, the Lebesgue constants of the single-parameter sets in both lower dimensions and in this work are reasonably small, thereby bringing the necessity of further optimization into question.

Finally, in the interest of completeness, the interpolation points found in this work can be tested independently of their Lebesgue constants. For example, one could perform a four-dimensional collocation projection of a transcendental function and compute the $L^2$ and/or $L^{\infty}$ error of the interpolation. We anticipate that this type of practical analysis will be investigated in future studies. 

\section*{Declaration of Competing Interests}

The authors declare that they have no known competing financial interests or personal relationships that could have appeared to influence the work reported in this paper.

\section*{Funding}

This research received funding from the United States Naval Research Laboratory (NRL) under grant number N00173-22-2-C008. In turn, the NRL grant itself was funded by Steven Martens, Program Officer for the Power, Propulsion and Thermal Management Program, Code 35, in the United States Office of Naval Research.


{\scriptsize\bibliography{technical-refs}}



\appendix

\section{Interpolation Point Tables} \label{point_tabs}

\scriptsize
    
\subsection{Optimized Interpolation Points, p = 1}~
\begin{longtable}{|r|r|r|r|r|}
    \hline & & & &\\[-2.0ex]
    $\lambda^1$ & $\lambda^2$ & $\lambda^3$ & $\lambda^4$ & $\lambda^5$ \\
    \hline & & & &\\[-2.0ex]
    1.000000000000000 & 0 & 0 & 0 & 0 \\
    0 & 1.000000000000000 & 0 & 0 & 0 \\
    0 & 0 & 1.000000000000000 & 0 & 0 \\
    0 & 0 & 0 & 1.000000000000000 & 0 \\
    0 & 0 & 0 & 0 & 1.000000000000000 \\
    \hline
\end{longtable}

\vspace{5mm}

\subsection{Optimized Interpolation Points, p = 2}~\begin{longtable}{|r|r|r|r|r|}
    \hline & & & &\\[-2.0ex]
    $\lambda^1$ & $\lambda^2$ & $\lambda^3$ & $\lambda^4$ & $\lambda^5$ \\
    \hline & & & &\\[-2.0ex]
    1.000000000000000                  & 0                  & 0                  & 0                  & 0 \\
    0.500000000000000  & 0.500000000000000                  & 0                  & 0                  & 0 \\
                    0  & 1.000000000000000                  & 0                  & 0                  & 0 \\ 
    0.500000000000000                  & 0  & 0.500000000000000                  & 0                  & 0 \\
                    0  & 0.500000000000000  & 0.500000000000000                  & 0                  & 0 \\
                    0                  & 0  & 1.000000000000000                  & 0                  & 0 \\
    0.500000000000000                  & 0                  & 0  & 0.500000000000000                  & 0 \\
                    0  & 0.500000000000000                  & 0  & 0.500000000000000                  & 0 \\
                    0                  & 0  & 0.500000000000000  & 0.500000000000000                  & 0 \\
                    0                  & 0                  & 0  & 1.000000000000000                  & 0 \\
    0.500000000000000                  & 0                  & 0                  & 0  & 0.500000000000000 \\
                    0  & 0.500000000000000                  & 0                  & 0  & 0.500000000000000 \\
                    0                  & 0  & 0.500000000000000                  & 0  & 0.500000000000000 \\
                    0                  & 0                  & 0  & 0.500000000000000  & 0.500000000000000 \\
                    0                  & 0                  & 0                  & 0  & 1.000000000000000 \\
    \hline
\end{longtable}

\vspace{5mm}

\subsection{Optimized Interpolation Points, p = 3}~\begin{longtable}{|r|r|r|r|r|}
    \hline & & & &\\[-2.0ex]
    $\lambda^1$ & $\lambda^2$ & $\lambda^3$ & $\lambda^4$ & $\lambda^5$ \\
    \hline & & & &\\[-2.0ex]
    1.000000000000000                  & 0                  & 0                  & 0                  & 0 \\
    0.723606797749979  & 0.276393202250021                  & 0                  & 0                  & 0 \\
    0.276393202250021  & 0.723606797749979                  & 0                  & 0                  & 0 \\
                    0  & 1.000000000000000                  & 0                  & 0                  & 0 \\
    0.723606797749979                  & 0  & 0.276393202250021                  & 0                  & 0 \\
\hline
\end{longtable}

\vspace{5mm}

\subsection{Optimized Interpolation Points, p = 4}~\begin{longtable}{|r|r|r|r|r|}
    \hline & & & &\\[-2.0ex] 
    $\lambda^1$ & $\lambda^2$ & $\lambda^3$ & $\lambda^4$ & $\lambda^5$ \\
    \hline & & & &\\[-2.0ex] 
    1.000000000000000                  & 0                  & 0                  & 0                  & 0 \\
    0.827326835353989  & 0.172673164646011                  & 0                  & 0                  & 0 \\ 
    0.500000000000000  & 0.500000000000000                  & 0                  & 0                  & 0 \\
    0.172673164646011  & 0.827326835353989                  & 0                  & 0                  & 0 \\
                    0  & 1.000000000000000                  & 0                  & 0                  & 0 \\
\hline
\end{longtable}

\vspace{5mm}

\subsection{Optimized Interpolation Points, p = 5}~\begin{longtable}{|r|r|r|r|r|}
    \hline & & & &\\[-2.0ex] 
    $\lambda^1$ & $\lambda^2$ & $\lambda^3$ & $\lambda^4$ & $\lambda^5$ \\
    \hline & & & &\\[-2.0ex] 
    1.000000000000000                  & 0                  & 0                  & 0                  & 0 \\
    0.882527661964732  & 0.117472338035268                  & 0                  & 0                  & 0 \\
    0.642615758240323  & 0.357384241759677                  & 0                  & 0                  & 0 \\
    0.357384241759677  & 0.642615758240323                  & 0                  & 0                  & 0 \\
    0.117472338035268  & 0.882527661964732                  & 0                  & 0                  & 0 \\
\hline
\end{longtable} 

\vspace{5mm}

\subsection{Optimized Interpolation Points, p = 6}~\begin{longtable}{|r|r|r|r|r|}
    \hline & & & &\\[-2.0ex] 
    $\lambda^1$ & $\lambda^2$ & $\lambda^3$ & $\lambda^4$ & $\lambda^5$ \\
    \hline & & & &\\[-2.0ex] 
    1.000000000000000                  & 0                  & 0                  & 0                  & 0 \\
    0.915111948139283  & 0.084888051860717                  & 0                  & 0                  & 0 \\ 
    0.734424396735357  & 0.265575603264643                  & 0                  & 0                  & 0 \\
    0.500000000000000  & 0.500000000000000                  & 0                  & 0                  & 0 \\ 
    0.265575603264643  & 0.734424396735357                  & 0                  & 0                  & 0 \\
\hline 
\end{longtable}

\end{document}